%
%
%
%
\input amstex
\advance\vsize-0.5cm\voffset=-0.5cm\advance\hsize1cm\hoffset0cm
\magnification=\magstep1
\documentstyle{amsppt}
\NoBlackBoxes
\topmatter
\title ON THE SPLITTING PROBLEM FOR MANIFOLD PAIRS WITH BOUNDARIES
\endtitle
\author Matija Cencelj -- Yuri V. Muranov -- Du\v san Repov\v s
\endauthor

\abstract The problem of splitting a homotopy equivalence along a submanifold
is
closely related
to the
surgery exact sequence and to the problem
of surgery of manifold pairs.
In  classical surgery theory there exist two approaches to surgery
in the category of manifolds with boundaries.
In the 
$ rel \ \partial$ 
case the surgery on a manifold
pair is
considered with the given fixed manifold structure
 on the boundary.
In the relative case the  surgery
on the manifold  with boundary is considered  without fixing maps on the
boundary. Consider a normal map to a manifold pair
$(Y, \partial Y)\subset (X, \partial X)$ with boundary
which is a simple homotopy equivalence
on
the boundary $\partial X$. This map defines a mixed structure 
on the manifold with the boundary in the sense of Wall.
We introduce and
study groups of obstructions to splitting of such 
mixed structures 
along submanifold with boundary $(Y, \partial Y)$.
We describe relations of these groups to classical surgery
and splitting  obstruction groups.
We also consider several geometric examples.
\endabstract
\endtopmatter

\document

\bigskip

\subhead 1. Introduction
\endsubhead
\smallskip

Let $(X^n, \partial X)$ be a  compact 
 topological $n$--manifold
with boundary. The set
$\Cal S^{CAT}(X,\partial X)$ of $CAT$-manifold structures
($CAT = TOP, PL, DIFF$) on
$(X, \partial X)$ consists of the classes of concordance
of simple homotopy equivalences of pairs
$f\: (M, \partial M)\to  (X, \partial X)$,
where $(M, \partial M)$ is a compact $CAT$-manifold pair of dimension $n$
with boundary (see \cite{7}, \cite{8}, and \cite{11}).
If $\partial X$ already has a $CAT$@-manifold  structure
then the set of  manifold structures on $X$ which
are fixed on the boundary is denoted by
$\Cal S^{CAT}_{\partial}(X,\partial X)$.

Let  $\Cal T^{CAT}(X, \partial X)$ be the set of classes
of normal bordisms of  normal maps
to the pair $(X, \partial X)$ and
$\Cal T^{CAT}_{\partial}(X,\partial X)$ the set of $rel \ \partial$
classes of normal bordisms of
normal maps (see \cite{7}, \cite{8} and \cite{11}).

Let $Y\subset X$ be a
submanifold  of a closed manifold $X^n$
of codimension
$q$.
Given a normal map $(f,b):M^n\to X^n$,
there is a problem
of finding
a simple homotopy equivalence $g:M_1\to X$
in the class of normal bordism $[(f,b)]$,
which is transversal to $Y$ and such that
$N=g^{-1}(Y)$ is a submanifold of $M_1$ 
and the restrictions
$$
g|_{N}:N\to  Y, \
g|_{M_1\setminus N}:M_1\setminus N\to X\setminus Y
\tag 1.1
$$
are
simple homotopy equivalences. The obstruction group
for doing such surgery is denoted by
$LP_{n-q}(F)$ (see \cite{11} and \cite{8}),
where
$$
F=\left(\matrix {\pi}_{1}(S(\xi)) & \longrightarrow &
{\pi}_{1}(X\setminus Y) \cr
\downarrow & \  & \downarrow \cr
{\pi}_{1}(D(\xi)) & \longrightarrow & {\pi}_{1}(X)
\endmatrix \right)
\tag 1.2
$$
is a pushout square  of fundamental groups with
orientations and $S(\xi)$ is the boundary of 
a
tubular neighborhood $D(\xi)$
of $Y$ in $X$.

If $f:M\to X$ is a simple homotopy equivalence
then the
obstruction to finding a map 
in the
homotopy class of the map $f$
with properties (1.1), 
which  is transversal to $Y$,
lies in the splitting  obstruction group $LS_{n-q}(F)$
(see \cite{11, \S 11} and \cite{8, \S 7.2}).

Let
$
(Y,\partial Y)\subset (X,\partial X)
$
be a codimension $q$ manifold pair with boundary.
In the 
$rel \ \partial$ case
the set  $\Cal T^{CAT}_{\partial}(X,\partial X)$
of tangent structures consists of classes of concordance rel 
boundary of normal maps
$$
(f,\partial f): (M,\partial M)\to (X, \partial X)
$$
with a 
fixed $CAT$-isomorphism 
$$
\partial f: \partial M\to \partial X
$$
which is already split on the boundary.
We have
a map
$$
\Cal T^{CAT}_{\partial}(X,\partial X)\to LP_{n-q}(F)
\tag 1.3
$$
which is given by
mapping the
obstruction to surgery to the normal map of manifold
pairs rel boundary.

In a similar way we
 have a map
$$
\Cal S^{CAT}_{\partial}(X,\partial X)\to LS_{n-q}(F)
\tag 1.4
$$
to the splitting obstruction group.

It follows from \cite{11, \S 11, page 136} (see also \cite{12})
that for the relative case we  have  maps
$$
\Cal T^{CAT}(X,\partial X) \to
LP_{n-q}(F_{\partial}\to F)
\tag 1.5
$$
similarly to (1.3) and (1.4) and
$$
  \Cal S^{CAT}(X,\partial X)\to LS_{n-q}(F_{\partial}\to F)
  \tag 1.6
$$
to the
relative obstruction groups
where $F_{\partial}$ is a pushout
square for a splitting problem of
the pair $\partial Y\subset \partial X$.

In accordance with Wall \cite{11, \S 10, p.~116} (see also \cite{3}),
it is possible to
introduce
a mixed type of structures on a manifold with boundary  $(X, \partial X)$.
Consider  a normal map
$$
(f,\partial f): (M,\partial M)\to (X, \partial X)
$$
for which the map
$$
\partial f: \partial M\to \partial X
$$
is a simple homotopy equivalence.
Two such maps  are concordant if they are normally bordant
by a bordism for which a restriction to  bordism between the boundaries
is an equivalence
in $S^{CAT}(\partial X)$.
Denote the set of concordance classes
by $\Cal T\Cal S^{CAT}(X, \partial X)$. The elements of
$\Cal T\Cal S^{CAT}(X, \partial X)$ are called mixed structures
on $(X,\partial X)$.

In the present paper we shall work in $TOP$@-category
and simple surgery obstruction groups (see \cite{7} and  \cite{8}).
%
We think of all surgery and splitting obstruction groups as 
decorated by an "s"
although we do not write it.

We introduce groups
$LPS_*(F_{\partial}\to F)$ and define a map
$$
\psi :\Cal T\Cal S(X, \partial X)\to LPS_*(F_{\partial}\to F)
$$
which gives an obstruction to finding a map
in the class of concordance
 which is split along the submanifold pair
$
(Y,\partial Y)\subset (X,\partial X).
$

We study properties of the introduced groups and their relations
to surgery and splitting obstruction groups. The main results are
given by braids of exact sequences. Then we consider several geometric
examples in which we compute the introduced groups
and natural maps.

In Section 2 we give explicit definitions of several structure sets
and recall the necessary technical results about the
algebraic
surgery exact sequences of Ranicki and the
surgery $L$@-spectrum.

In Section 3  we recall main properties of
splitting obstruction groups and introduce $LPS_*$@-groups.
These groups are  realized
as homotopy groups of a spectrum. We describe algebraic properties
of these groups and  relations of these
groups to surgery and splitting obstruction groups and
to surgery exact sequence.

In Section 4 we consider  geometric examples
in which we compute $LPS_*$@-groups and natural maps
between introduced groups and classical obstruction groups
which arise naturally in the considered problem.
\bigskip

\subhead 2. Structure sets and surgery exact sequence
\endsubhead
\smallskip

For definitions of structure sets we shall follow  Ranicki
\cite{8}.  Let $X^n$ be a closed topological manifold.
A $t$@-triangulation of $X$ is a topological normal map
(see \cite{8, \S 1.2})
$$
(f,b): M\to X,
$$
where $M$ is a closed $n$@-dimensional topological manifold.
Two $t$@-triangulations
$$
(f_i,b_i): M_i\to X, \  i=0,1
$$
are concordant  \cite{8, \S 7.1}
  if there exists a topological normal
map of triads
$$
((g, c);(f_0,b_0), (f_1,b_1): (W; M_0, M_1)\to (X\times I; X\times\{0\},
X\times\{1\})
$$
where $I=[0,1]$ and $W$ is a compact $(n+1)$@-dimensional manifold
with boundary $\partial W=M_0\cup M_1$.
The set of concordance classes of $t$@-triangulations of $X$ is
denoted by $\Cal T^{TOP}(X)$. Note that we shall consider the case
of a manifold $X$ and hence the set  $\Cal T^{TOP}(X)$ will be
nonempty.

An $s$@-triangulation of a closed topological manifold $X^n$
is a simple homotopy equivalence
$
f:M\to X,
$
where $M$ is a closed topological $n$@-dimensional manifold.

Two $s$@-triangulations
$$
(f_i,b_i): M_i\to X, \  i=0,1
$$
are concordant  \cite{8, \S 7.1}
  if there exists a simple homotopy equivalence  of triads
$$
(g;f_0, f_1): (W; M_0, M_1)\to (X\times I; X\times\{0\},
X\times\{1\})
$$
where $W$ is a compact $(n+1)$@-dimensional manifold
with the boundary $\partial W=M_0\cup M_1$.
The set of concordance classes of $s$@-triangulations of $X$ is
denoted by $\Cal S^{TOP}(X)$. This set is called
the {\it topological manifold structure set}. As before,
the set  $\Cal S^{TOP}(X)$ will be
nonempty. These sets fit in the surgery exact sequence (see \cite{8, \S 7}
and \cite{11})
$$
\cdots L_{n+1}(\pi_1(X))\to \Cal S^{TOP}(X)\to \Cal T^{TOP}(X)\to
L_n(\pi_1(X))
\tag 2.1
$$
where $L_*(\pi_1(X))$ are surgery obstruction groups.

Now consider the case of a compact $n$@-dimensional manifold $X$  with
the boundary $\partial X$.
First, we  consider the case
of structures which are fixed on the boundary. This is
the $rel \ \partial$ case.
A $t_{\partial}$@-triangulation of $(X, \partial X)$
is a topological normal map of pairs (see \cite{8, \S 7.1})
$$
((f,b),(\partial f, \partial b)):
(M, \partial M)\to (X, \partial X)
$$
with a homeomorphism $\partial f:\partial M\to \partial X$.
Two $t_{\partial}$@-triangulations
$$
((f_i,b_i),(\partial f_i, \partial b_i)):
(M_i, \partial M_i)\to (X, \partial X), \ i=0,1
$$
are concordant
if there exists a topological normal map
$$
\matrix
((h,d);(g,c),(f_0, b_0), (f_1, b_1)):&&\\
(W; V, M_0, M_1)&\longrightarrow& (X\times I;
\partial X\times I, X\times\{0\},
 X\times\{1\})
\endmatrix
$$
with
$$
V=\partial M_0\times I, \
\partial V=\partial M_0\cup \partial M_1
$$
and
$$
(g,c)=\partial f_0 \times I:V\to \partial X\times I.
$$
The set of concordance classes is denoted
by $\Cal T^{CAT}_{\partial}(X,\partial X)$
(see \cite{11, \S 10} and \cite{8, \S 7.1}).

An $s_{\partial}$@-triangulation of $(X, \partial X)$
is a simple homotopy equivalence  of pairs (see \cite{8, \S 7.1})
$$
(f,\partial f):(M, \partial M)\to (X, \partial X)
$$
with a homeomorphism $\partial f:\partial M\to \partial X$.
Two $s_{\partial}$@-triangulations
$$
(f_i,\partial f_i):
(M_i, \partial M_i)\to (X, \partial X), \ i=0,1
$$
are concordant
if there exists a simple homotopy equivalence of $4$@-ads
$$
(h;g, f_0, f_1):
(W; V, M_0, M_1)\to (X\times I; \partial X\times I,  X\times\{0\},
 X\times\{1\})
$$
with
$$
V=\partial M_0\times I, \
\partial V=\partial M_0\cup \partial M_1
$$
and
$$
g=\partial f_0 \times I:V\to \partial X\times I.
$$
The set of concordance classes is denoted
by $\Cal S^{CAT}_{\partial}(X,\partial X)$
(see \cite{11, \S 10} and \cite{8, \S 7.1}).

These sets fit in the  surgery exact  sequence 
(see \cite{11, \S 10} and \cite{8, \S 7})
$$
\cdots \to   L_{n+1}(\pi_1(X)) \rightarrow
\Cal S^{TOP}_{\partial}(X, \partial X) \to
\Cal T^{TOP}_{\partial}(X, \partial X)
\to L_n(\pi_1(X)).
\tag 2.2
$$

Now consider the relative case
of structures on a manifold with boundary.
A $t$@-triangulation of $(X, \partial X)$
is a topological normal map of pairs (see \cite{8, \S 7.1})
$$
((f,b),(\partial f, \partial b)):
(M, \partial M)\to (X, \partial X).
$$
Two $t$@-triangulations
$$
((f_i,b_i),(\partial f_i, \partial b_i)):
(M_i, \partial M_i)\to (X, \partial X), \ i=0,1
$$
are concordant
if there exists a topological normal map of $4$@-ads
$$
((h,d);(g,c),(f_0, b_0), (f_1, b_1)):
(W; V, M_0, M_1)\to
(X\times I; \partial X\times I,  X\times\{0\},
 X\times\{1\})
$$
with
$$
\partial V=\partial M_0\cup \partial M_1.
$$
The set of concordance classes is denoted
by $\Cal T^{CAT}(X,\partial X)$
(see \cite{11, \S 10} and \cite{8, \S 7.1}).

An $s$@-triangulation of $(X, \partial X)$
is a simple homotopy equivalence  of pairs (see \cite{8, \S 7.1})
$$
(f,\partial f):(M, \partial M)\to (X, \partial X).
$$
Two $s$@-triangulations
$$
(f_i,\partial f_i):
(M_i, \partial M_i)\to (X, \partial X), \ i=0,1
$$
are concordant
if there exists a simple homotopy equivalence of $4$@-ads
$$
((h;g, f_0, f_1):
(W; V, M_0, M_1)\to (X\times I; \partial X\times I,  X\times\{0\},
 X\times\{1\})
$$
with
$$
\partial V=\partial M_0\cup \partial M_1.
$$
The set of concordance classes is denoted
by $\Cal S^{CAT}(X,\partial X)$
(see \cite{11, \S 10} and \cite{8, \S 7.1}).

These sets fit in the  surgery exact  sequence
(see \cite{11, \S 10} and \cite{8, \S 7})
$$
\cdots \rightarrow
\Cal S^{TOP}(X, \partial X) \to
\Cal T^{TOP}(X, \partial X)
\to L_n\left(\pi_1(\partial X)\to \pi_1(X)\right).
\tag 2.3
$$

Now we define mixed structures on a manifold with boundary
(see \cite{11, page 116} and \cite{3}).
A $t_s$@-triangulation of $(X, \partial X)$
is a topological normal map of pairs
$$
((f,b),(\partial f, \partial b)):
(M, \partial M)\to (X, \partial X)
$$
such that $\partial f:\partial M\to \partial X$ is
an $s$@-triangulation.
Two $t_s$@-triangulations
$$
((f_i,b_i),(\partial f_i, \partial b_i)):
(M_i, \partial M_i)\to (X, \partial X), \ i=0,1
$$
are concordant
if there exists a topological normal map
$$
\matrix
((h,d);(g,c),(f_0, b_0), (f_1, b_1)):&&\\
(W; V, M_0, M_1)&\longrightarrow& (X\times I;
\partial X\times I,  X\times\{0\},
 X\times\{1\})
\endmatrix
$$
with
$$
\partial V=\partial M_0\cup \partial M_1
$$
and
$
g:V\to \partial X\times I
$
is a concordance of $s$@-triangulations $\partial f_0$ and $\partial f_1$.
The set of concordance classes is denoted
by $\Cal T\Cal S^{TOP}(X,\partial X)$
(see \cite{11, page 116}).

It follows  from definitions (see also \cite{3}) that  the following natural
forgetful maps
$$
\Cal T\Cal S^{TOP}(X, \partial X)
\to \Cal T^{TOP}(X, \partial X),
$$
$$
\Cal S^{TOP}(X, \partial X)\to
\Cal T\Cal S^{TOP}(X, \partial X),
\tag 2.4
$$
$$
 \Cal T\Cal S^{TOP}(X, \partial X) \to
\Cal S^{TOP}(\partial X)
$$
are well--defined.

The maps in (2.4) fit  in the following
exact sequences  (see \cite{11, page 116},
\cite{8, \S 7}, and \cite{3})
$$
\cdots \to L_n(\pi_1(\partial X))\to \Cal T\Cal S^{TOP}(X, \partial X)
\to \Cal T^{TOP}(X, \partial X) \to L_{n-1}(\pi_1(\partial X)),
\tag 2.5
$$
$$
\cdots \to L_{n+1}(\pi_1(X))\to \Cal S^{TOP}(X, \partial X)
\to \Cal T\Cal S^{TOP}(X, \partial X) \to L_{n}(\pi_1(X)),
\tag 2.6
$$
and
$$
\cdots \to
\Cal T_{\partial}^{TOP}(X, \partial X)\to
 \Cal T\Cal S^{TOP}(X, \partial X) \to
\Cal S^{TOP}(\partial X).
\tag 2.7
$$

We now recall the necessary results concerning the application
of homotopy category of spectra to surgery theory
(see \cite{1}, \cite{2}, \cite{4}, \cite{5}, \cite{6},
and \cite{7}).
A spectrum $\Bbb E$ is given by a collection of
$CW$-complexes $\{(E_n,*)\}$, $n \in \Bbb Z$, together with
cellular  maps
$\{\epsilon_n:SE_n\to E_{n+1}\}$,
where  $SE_n$ is the suspension of the space $E_n$ (see \cite{14}).
A spectrum $\Bbb E$ is an $\Omega$@-spectrum if the
adjoint maps $\epsilon_n^{\prime}:E_n\to \Omega E_{n+1}, \ n\in \Bbb Z$
are homotopy equivalences.

In the category of spectra the suspension functor $\Sigma$
and iterated functors $\Sigma^{k}$, $k\in \Bbb Z$
are well-defined (see \cite{10}).
For every spectrum  $\Bbb E$ we have  an isomorphism of homotopy
groups  $\pi_n(\Bbb E)=\pi_{n+k}(\Sigma^k \Bbb E)$.
Recall that  in the homotopy category of spectra the concepts of pull-back
and push-out squares are equivalent.

In accordance with \cite{7}, \cite{8}, and  \cite{11}
the surgery obstruction groups
$L_n(\pi)$  and such natural maps as induced by inclusion
and transfer are realized on the
spectrum level.
That is,
for every group $\pi$ with a homomorphism of
orientation $\omega: \pi\to \{\pm 1\}$ there exists an $\Omega$-spectrum
$\Bbb L(\pi, \omega)$ with homotopy groups
$$
\pi_n(\Bbb L(\pi, \omega)) = L_n(\pi, \omega).
$$
In what follows we shall not include homomorphism of orientation
in our notation and will assume that all groups are equipped with
such a homomorphism and all homomorphisms of groups preserve
orientation.
Any homomorphism of groups  $f : \pi \to G$
induces a cofibration of spectra
$$
\Bbb L(\pi) \to \Bbb L(G) \to \Bbb L(f)
\tag 2.8
$$
where $\Bbb L(f)$ is the spectrum for the
relative $L$@-groups
$$
L_n(f)=L_n(\pi\to G)=\pi_n(\Bbb L(f)).
$$
We have a similar situation  for the transfer map
(see for example, \cite{11} and \cite{12}).

Let  $X$ be a topological space.
An algebraic surgery  exact sequence of Ranicki
(see \cite{7} and \cite{8})
$$
\cdots\rightarrow L_{n+1}(\pi_1(X))\rightarrow \Cal  S_{n+1}(X)
\rightarrow
H_n(X, \bold L_{\bullet})\rightarrow L_{n}(\pi_1(X))\rightarrow
\cdots
\tag 2.9
$$
is defined.
Here  $\bold L_{\bullet}$ is  the 1-connected cover of the
surgery $\Omega$@-spectrum
$\Bbb L(1)$ with
$\{\bold L_{\bullet}\}_0 \simeq G/TOP$.
The algebraic surgery exact  sequence (2.9)
is the homotopy long exact sequence of the cofibration
$$
X_+\land \bold L_{\bullet}\rightarrow \Bbb  L(\pi_1(X)).
\tag 2.10
$$

By definition, we have $\Cal S_i(X)=\pi_i(\Bbb S(X))$ for the homotopy
cofiber $\Bbb S(X)$ of the map in (2.10).
For a closed $n$@-dimensional topological manifold $X$
we have
$$
\pi_{n+1}(\Bbb S(X)) =\Cal S_{n+1}(X)\cong \Cal S^{TOP}(X),\
\Cal T(X)\cong H_n(X; L_{\bullet}),
\tag 2.11
$$
and the surgery exact sequence (2.1) is isomorphic to the left part
of the algebraic surgery exact sequence (2.9)
(see \cite{8, Proposition 7.1.4}).

For the case of a
compact topological manifold  $X$ with boundary $\partial X$ the algebraic
surgery exact sequences for the relative case
and for the $rel \ \partial$ case
are contained 
in  the following commutative diagram (see \cite{7} and  \cite{8, \S 7})
$$
\matrix
\ \ \     \vdots        &     &\vdots          &   &
\vdots          &     & \vdots   &     \\
\ \ \     \downarrow        &     &\downarrow          &   &
\downarrow          &     & \downarrow   &     \\
\cdots L_{n+1}(\pi)& \to &\Cal S^{\partial}_{n+1}(X,\partial X)
&\to&
H_n(X; \bold L_{\bullet})& \to & L_n(\pi)&\cdots  \\
\ \ \    \downarrow        &     &\downarrow
&   & \downarrow          &     & \downarrow   &      \\
\cdots L_{n+1}^{rel}& \to&
\Cal S_{n+1}(X,\partial X) &\to
&H_n(X,\partial X;\bold L_{\bullet})&
\to &L_n^{rel}&\cdots  \\
\ \ \     \downarrow        &     &\downarrow          &   &
\downarrow          &     & \downarrow   &                \\
\cdots L_{n}(\rho)& \to &\Cal S_{n}(\partial X)
&\to&
H_{n-1}(\partial X; \bold L_{\bullet})& \to &
L_{n-1}(\rho)&\cdots\\
\ \ \     \downarrow        &     &\downarrow          &   &
\downarrow          &     & \downarrow   &     \\
\ \ \     \vdots        &     &\vdots          &   &
\vdots          &     & \vdots   &     \\
\endmatrix
\tag 2.12
$$
where $\pi=\pi_1(X)$,  $\rho=\pi_1(\partial X)$, and
$L_*^{rel}=L_*(\rho\to\pi)$.
Diagram (2.12)
is realized on the
spectrum level (see \cite{8}, \cite{1}, and \cite{3}).
We denote by $\Bbb S(X,\partial X)$ the homotopical cofiber 
of the map 
$$
(X/\partial X)_+\land\bold L_{\bullet}\to \Bbb L(\rho\to \pi),
$$
and by $\Bbb S^{\partial}(X,\partial X)$ the homotopical cofiber 
of the map 
$$
X_+\land\bold L_{\bullet}\to \Bbb L(\pi).
$$
We have 
$$
\pi_i(\Bbb S^{\partial}(X, \partial X))=
\Cal S^{\partial}_i(X, \partial X)
$$
and 
$$
\pi_i(\Bbb S(X, \partial X))= \Cal S_i(X, \partial X).
$$
For a topological manifold $X$
the left part of the upper row of diagram (2.12) is isomorphic to
the exact sequence (2.2).
The left part of middle row of diagram (2.12) is isomorphic to
exact sequence (2.3).

In particular, we have the
isomorphisms
$$
S^{TOP}(X, \partial X)\cong S_{n+1}(X, \partial X), \quad
S^{TOP}_{\partial}(X, \partial X)\cong  S^{\partial}_{n+1}(X, \partial X),
\tag 2.13
$$
and
$$
\Cal T^{TOP}(X, \partial X)\cong H_n(X, \partial X; \bold L_{\bullet}),
\quad \Cal T_{\partial}^{TOP}(X, \partial X)\cong H_n(X; \bold L_{\bullet}).
\tag 2.14
$$

Consider the
composition
$$
L_{n+1}(\pi_1(X))\to L_{n+1}(\pi_1(\partial X)\to \pi_1(X))
\to \Cal S^{TOP}(X, \partial X)
\tag 2.15
$$
where the first map lies in the relative exact sequence
of $L$-groups for the map $\pi_1(\partial X)\to \pi_1(X)$
and the second map lies in (2.3). It follows from (2.12)
that
the composition
(2.15) is realized by a map of spectra (see also \cite{3})
$$
\Bbb L(\pi_1(X)) \to \Bbb S(X, \partial X)
\tag 2.16
$$
and the cofiber of the map in (2.16) is denoted
by $\Bbb T\Bbb S(X,\partial X)$.
We shall denote
$$
\pi_n(\Bbb T\Bbb S(X,\partial X)) =\Cal T\Cal S_n(X,\partial X)
$$
and we get
an isomorphism \cite{3}
$$
\Cal T\Cal S_{n+1}(X,\partial X)\cong\Cal T\Cal S^{TOP}(X,\partial X).
$$
Note that in a similar way (see \cite{3}) it is possible to
describe the spectrum $\Bbb T\Bbb S(X,\partial X)$
as the homotopical cofiber of any of the following maps
$$
\Bbb S(\partial X)\to \Sigma(X_+\land\bold L_{\bullet})
\ \text{and} \
(X/\partial X)_+\land\bold L_{\bullet})\to\Sigma \Bbb L(\pi_1(\partial X)).
\tag 2.17
$$

\bigskip

\subhead 3. Splitting problem for a manifold with boundary
\endsubhead
\smallskip

Let $(X, Y, \xi)$ be a codimension $q (=1,2)$   manifold
pair  in the sense of Ranicki (see \cite{8, page 570}), i.e.
a locally flat closed submanifold  $Y$ is given with a normal block bundle
 $$
  \xi=\xi_{Y\subset X}: Y \to \widetilde{BTOP}(q)
 $$
for which we have a decomposition 
of the closed manifold
$$
 X =D(\xi)\cup_{S(\xi)}\overline{X\setminus D(\xi)}.
$$
where $D(\xi)$ is the total space of the normal block bundle
with the boundary $S(\xi)$.
In accordance with \cite{8, p. 570} the pair $(X,Y)$ has an underlying 
 structure
of  an
$(n, n-q)$@-dimensional $t$@-normal geometric Poincar\'e pair 
with the associated $(D^q, S^{q-1})$
fibration
$$
(D^q, S^{q-1})\to (D(\xi), S(\xi))\to Y.
\tag 3.1
$$

The fibration (3.1) provides transfer maps  on the
spectrum level
(see \cite{1}, \cite{6},
\cite{11}, and \cite{12})
$$
p^{\sharp}:\Bbb L(\pi_1(Y)) \to \Omega^q\Bbb
L\left(\pi_1(S(\xi))\to \pi_1(D(\xi))\right)
\tag 3.2
$$
and
$$
p^{\sharp}_1:\Bbb L(\pi_1(Y)) \to\Omega^{q-1}\Bbb L
\left(\pi_1(S(\xi))\right).
\tag 3.3
$$
Transfer  maps (3.2) and (3.3) fit in a homotopy commutative diagram
of spectra
$$
\matrix
\Bbb L(\pi_1(Y)) & \overset{p^{\sharp}}\to{\rightarrow}
&\Omega^q\Bbb
L(\pi_1(S(\xi))\to \pi_1(D(\xi)))&
\to & \Omega^q\Bbb
L(\pi_1(X\setminus Y)\to \pi_1(X))  \\
  & \ p^{\sharp}_1\searrow & \downarrow & & \downarrow \\
 & & \Omega^{q-1}\Bbb L(\pi_1(S(\xi)))&
 \to &\Omega^{q-1}\Bbb
L(\pi_1(X\setminus Y)),
\endmatrix
\tag 3.4
$$
where the horizontal
maps of the right square
are induced by the horizontal maps of $F$ and the vertical maps are
maps from cofibrations of spectra as in (2.8) for the vertical maps of
the square $F$.

The spectrum
$
\Bbb L\Bbb S (F)
$
for splitting obstruction groups of the manifold pair $Y\subset X$ and
the spectrum
$
\Bbb L\Bbb P(F)
$
for surgery obstruction groups of the manifold pair
fit in the homotopy commutative diagram of spectra
$$
\CD
\Omega\Bbb L(\pi_1(Y)) @>>>
\Omega^{q+1}\Bbb L(\pi_1(X\setminus Y)\to \pi_1(X))
@>>> \Bbb L\Bbb  S(F) \\
@V = VV @VVV @VVV   \\
\Omega\Bbb L(\pi_1(Y)) @>>> \Omega^q\Bbb L(\pi_1(X\setminus Y))
@>>> \Bbb  L \Bbb P(F)
\endCD
\tag 3.5
$$
where the left horizontal maps are compositions from
diagram (3.4) and the right square is the
pullback (see \cite{1}, \cite{10}, and \cite{11}).
In particular,  we have the
isomorphisms
$$
\pi_n(\Bbb L\Bbb  S(F))\cong LS_n(F), \
\pi_n(\Bbb L\Bbb  P(F))\cong LP_n(F).
$$

A topological normal map  \cite{8, \S 7.2}
$$
((f,b), (g,c)):(M,N)\to (X,Y)
$$
to the manifold pair $(X, Y, \xi)$ is
represented by a normal map
$(f,b)$ to the manifold $X$
which is transversal to $Y$ with $N=f^{-1}(Y)$, and
$(M,N)$ is a topological manifold pair with the normal
block bundle 
$$
\nu:N\overset{f|_N}\to{\to} Y \overset{\xi}\to{\to}
\widetilde{BTOP}(q).
$$
Additionally, the following conditions are satisfied:

(i) the restriction
$$
(f,b)|_N =(g,c) :N\to Y
$$
is a normal map;

(ii) the restriction
$$
(f,b)|_P =(h,d) :(P, S(\nu))\to (Z,S(\xi))
$$
is a normal map to the pair $(Z,S(\xi))$,
where
$$
P=\overline{M\setminus D(\nu)}, \ \ Z=\overline{X\setminus
D(\xi)};
$$

(iii) the restriction
$$
(h,d)|_{S(\nu)}: S(\nu)\to S(\xi)
$$
coincides with the induced map
$$
(g,c)^!: S(\nu)\to S(\xi),
$$
and $(f,b)=(g,c)^!\cup(h,d)$.

The normal maps to $(X,Y,\xi)$ are called $t$@-triangulations
of the manifold pair $(X,Y)$
and the set of concordance classes of $t$@-triangulations
of the pair $(X, Y, \xi)$ coincides with the set of $t$@-triangulations
of the manifold $X$ \cite{8, Proposition 7.2.3}.

 An $s$@-triangulation of a manifold pair $(X,Y, \xi)$  in
topological category
 \cite{8, p. 571}
is a $t$@-triangulation
 of this pair for which the maps
 $$
 f:M\to X,\ g:N\to Y, \
 \text{and} \
 (P, S(\nu))\to (Z, S(\xi))
\tag 3.6
 $$
are simple homotopy equivalences ($s$@-triangulations).
The set of concordance classes of $s$@-triangulations is denoted
by $\Cal S^{TOP}(X, Y, \xi)$ (see \cite{8, page 571}).
Natural forgetful maps
$$
\Cal S^{TOP}(X, Y, \xi)\to \Cal S^{TOP}(X) \ \text{and} \
\Cal S^{TOP}(X, Y, \xi)\to \Cal T^{TOP}(X)
\tag 3.7
 $$
are well-defined (see \cite{8, \S 7.2}).
We have also the maps of taking obstruction (see \cite{8, page 572})
$$
\Cal S^{TOP}(X)\to LS_{n-q}(F) \ \text{and} \
\Cal T^{TOP}(X)\to LP_{n-q}(F).
\tag 3.8
$$
The maps in (3.7) and (3.8) are realized on the level of spectra 
(see \cite{1}, \cite{8 \S 7.2}, and  \cite{11}).
We shall denote by  $\Bbb S(X, Y, \xi)$ the homotopy cofiber 
of the map 
$$
X_+\land \bold L_{\bullet}\rightarrow \Sigma^q\Bbb  L\Bbb P(F)
$$
and by $\Cal S_i(X, Y, \xi)=\pi_i(\Bbb S(X, Y, \xi))$ 
its homotopy groups. 
We have an isomorphism 
$$
\Cal S_{n+1}(X, Y, \xi)\cong  S^{TOP}(X, Y, \xi).
\tag 3.9
$$
The maps from (3.7) and (3.8) fit in several diagrams
of exact sequences which are given in \cite{8, Proposition 7.2.6}.
The diagram
 $$
\matrix \rightarrow & L_{n+1}(\pi_1(X)) & \rightarrow &
LS_{n-q}(F) &
\rightarrow &  \Cal S_{n}(X,Y, \xi)& \rightarrow \cr
\ &  \nearrow \ \ \ \ \ \ \ \ \searrow & \ &  \nearrow \ \ \ \ \ \ \
\
\searrow
& \  & \nearrow \ \ \ \ \ \ \ \  \searrow & \ \cr
\ & \ & \Cal S_{n+1}(X)& \ &
LP_{n-q}(F) & \ & \ \cr
\ &  \searrow \ \ \ \ \ \ \ \ \nearrow & \ &  \searrow \ \ \ \ \ \ \
\
\nearrow
& \  & \searrow \ \ \ \ \ \ \ \  \nearrow & \ \cr
\rightarrow & \Cal S_{n+1}(X,Y, \xi ) & \longrightarrow &
H_{n}(X; \bold L_{\bullet}) &
\longrightarrow & L_n(\pi_1(X)) & \rightarrow
\endmatrix,
\tag 3.10
$$
from \cite{8, Proposition 7.2.6} is realized on the spectrum level
with the left part ($i\geq n$) which is isomorphic to a
geometrically
defined diagram  (see \cite{8, page 582}) containing structure sets
$\Cal S^{TOP}(X), \Cal S^{TOP}(X,Y, \xi), \Cal T^{TOP}(X),$
in accordance with  isomorphisms (2.11) and (3.9).
Note here that the geometric version of  diagram (3.10)
also 
contains the maps from (3.7) and (3.8).

Let
$$
(Y,\partial Y)\subset (X,\partial X)
\tag 3.11
$$
be a codimension $q$ manifold pair with boundary.
A manifold pair (3.11)
with boundaries defines a  pair of closed manifolds
$\partial Y\subset \partial X$ with
 a pushout square
$$
F_{\partial}=\left(\matrix {\pi}_{1}(S(\partial\xi)) & \longrightarrow &
{\pi}_{1}(\partial X\setminus \partial Y) \cr
\downarrow & \  & \downarrow \cr
{\pi}_{1}(\partial Y) & \longrightarrow & {\pi}_{1}(\partial X)
\endmatrix \right)
\tag 3.12
$$
of fundamental groups for the splitting problem.
A natural inclusion $\delta: \partial X \to X$ induces
a map of $\Delta:F_{\partial}\to F$ of squares of fundamental groups.

In the $rel \ \partial$@-case we
consider  $t$@-triangulations
$$
(f,\partial f): (M,\partial M)\to (X, \partial X)
\tag 3.13
$$
which are split on the boundary along $\partial Y$.
The classes of concordance relative to the boundary
of such maps give the set $\Cal T^{TOP}_{\partial}(X, \partial X)$
(see \cite{8, \S 7.2}) and the map
$$
\Cal T^{TOP}_{\partial}(X, \partial X)\to LP_{n-q}(F)
\tag 3.14
$$
defines a $rel\ \partial$ codimension $q$ splitting obstruction
along $Y\subset X$ (see \cite{8, \S 7.2}).

In a similar way  (see \cite{8, \S 7.2}) we
can consider an $s$@-triangulation
of pairs (3.12) which is split along the boundary.
The set of concordance $rel\ \partial$ classes is
$\Cal S^{TOP}_{\partial}(X, \partial X)$ and
a $rel \ \partial$ codimension $q$ splitting obstruction gives
a map
$$
\Cal S^{TOP}_{\partial}(X, \partial X)\to LS_{n-q}(F).
\tag 3.15
$$

As in the case of closed manifolds
denote by $\Cal S^{TOP}_{\partial}(X,Y, \xi)$ the set of
classes of concordance $rel\ \partial$ maps which are
split along $Y\subset X$.

The algebraic version of surgery  exact sequence (2.2)
and algebraic versions of the maps (3.14) and (3.15)
fit in the commutative braid of exact sequences
 $$
\matrix \rightarrow & L_{n+1}(\pi_1(X)) & \rightarrow &
LS_{n-q}(F) &
\rightarrow &  \Cal S_{n}^{\partial}(X,Y, \xi)& \rightarrow \cr
\ &  \nearrow \ \ \ \ \ \ \ \ \searrow & \ &  \nearrow \ \ \ \ \ \ \
\
\searrow
& \  & \nearrow \ \ \ \ \ \ \ \  \searrow & \ \cr
\ & \ & \Cal S_{n+1}^{\partial}(X, \partial X)& \ &
LP_{n-q}(F) & \ & \ \cr
\ &  \searrow \ \ \ \ \ \ \ \ \nearrow & \ &  \searrow \ \ \ \ \ \ \
\
\nearrow
& \  & \searrow \ \ \ \ \ \ \ \  \nearrow & \ \cr
\rightarrow & \Cal S_{n+1}^{\partial}(X,Y, \xi ) & \longrightarrow &
H_{n}(X; \bold L_{\bullet}) &
\longrightarrow & L_n(\pi_1(X)) & \rightarrow
\endmatrix
\tag 3.16
$$
Diagram (3.16) is realized on the level of spectra and for $\partial
X=\emptyset$
coincides with the diagram (3.10).
The left part ($i\geq n$) of diagram (3.16) is isomorphic to geometrically
defined diagram similarly to diagram (3.10).
In particular,
$$
\Cal S_i^{\partial}(X,Y, \xi)=
\pi_i(\Bbb S^{\partial}(X, Y, \xi)),
\text{and} \
\Cal S^{\partial}_{n+1}(X, Y, \xi)\cong \Cal S^{TOP}_{\partial}(X,Y, \xi).
\tag 3.17
$$

Denote by $LS_*(\Delta)=LS_*(F_{\partial} \to F)$
and $LP_*(\Delta)=LP_*(F_{\partial} \to F)$ the relative groups for
the map of squares $\Delta : F_{\partial} \to F$ which is induced
by the natural inclusion $\delta:\partial X \to X$.
It follows from functoriality of diagram (3.5) that these relative groups
are realized on the level of spectra. We have
  cofibrations of $\Omega$@-spectra
$$
\Bbb L\Bbb S(F_{\partial})\to\Bbb L\Bbb S(F)\to
\Bbb L\Bbb S(\Delta)
\tag 3.18
$$
and
$$
\Bbb L\Bbb P(F_{\partial})\to\Bbb L\Bbb P(F)\to
\Bbb L\Bbb P(\Delta)
\tag 3.19
$$
where
$$
\pi_n(\Bbb L\Bbb S(\Delta))\cong
 LS_n(F_{\partial}\to F) \ \text{and} \
 \pi_n(\Bbb L\Bbb P(\Delta))\cong
 LP_n(F_{\partial}\to F).
$$
These groups
fit in the commutative diagram of exact sequences
$$
\matrix
\ \ \     \vdots        &     &\vdots          &   &
\vdots          &     & \vdots   &     \\
\ \ \     \downarrow        &     &\downarrow          &   &
\downarrow          &     & \downarrow   &     \\
\cdots L_{n+1}(\rho)& \to &LS_{n-q}(F_{\partial})
&\to&
LP_{n-q}(F_{\partial})& \to &
L_{n}(\rho)&\cdots\\
\ \ \     \downarrow        &     &\downarrow          &   &
\downarrow          &     & \downarrow   &     \\
\cdots L_{n+1}(\pi)& \to &LS_{n-q}(F)
&\to&
LP_{n-q}(F)& \to & L_n(\pi)&\cdots  \\
\ \ \    \downarrow        &     &\downarrow
&   & \downarrow          &     & \downarrow   &      \\
\cdots L_{n+1}(\rho\to\pi)& \to&
LS_{n-q}(\Delta) &\to
&LP_{n-q}(\Delta)&
\to &L_n(\rho\to\pi)&\cdots  \\
\ \ \     \downarrow        &     &\downarrow          &   &
\downarrow          &     & \downarrow   &                \\
\ \ \     \vdots        &     &\vdots          &   &
\vdots          &     & \vdots   &     \\
\endmatrix
\tag 3.20
$$
where $\pi=\pi_1(X)$ and  $\rho=\pi_1(\partial X)$.
Diagram (3.20)
is realized on the  level of spectra  and the two middle columns are homotopy long exact
sequences of
cofibrations (3.18) and (3.19).

Now consider  relative structure groups for
a codimension $q$ manifold pair with boundaries (3.11).
We  have a normal block bundle 
$(\xi, \partial \xi)$ over the pair $(Y, \partial Y)$
and a decomposition
$$
(X, \partial X)= (D(\xi)\cup_{S(\xi)}Z,
D(\partial\xi)\cup_{S(\partial\xi)}\partial_+Z)
\tag 3.21
$$
where $(Z; \partial_+Z, S(\xi); S(\partial\xi))$ is a manifold
triad. Note here that $\partial_+Z =\overline{\partial X\setminus
D(\partial\xi)}$.

A topological
normal map (3.13)
of manifold pairs with boundaries provides
a normal block bundle
$(\nu, \partial \nu)$ over the pair $(N, \partial N),$
where (see \cite{8, p. 570})
$$
(N, \partial N)=(f^{-1}(Y),(\partial f)^{-1}(\partial Y)).
$$
 We have the
following decomposition
$$
(M, \partial M)= (D(\nu)\cup_{S(\nu)}P,
D(\partial\nu)\cup_{S(\partial\nu)}\partial_+P)
\tag 3.22
$$
where $(P; \partial_+P, S(\nu); S(\partial\nu))$ is a manifold
triad.

Let
$$
(f,\partial f):(M,\partial M)\to
(X, \partial X)
$$
be a normal map of a codimension $q$ pair with boundary
$(N, \partial N)\subset (M, \partial M)$  to
a codimension $q$ pair $(Y, \partial Y)\subset (X, \partial X)$.
It is an $s$@-triangulation
if the maps
$f:M\to X$
and
$\partial f : \partial M \to  \partial X$
are  $s$@-triangulations
of corresponding  codimension $q$ pairs.
We shall denote the set of concordance classes of $s$@-triangulations
of the codimension $q$ manifold pair (3.11) by
$$
\Cal S^{TOP}(X, Y;\partial) =
\Cal S^{TOP}(X,\partial X; Y, \partial Y; \xi, \partial(\xi)).
$$
The relative surgery theory (see \cite{8, \S 7.2}, \cite{10, \S 11},
and \cite{11}) guarantees that this structure set fits
in the following exact sequences.
$$
\cdots \to \Cal S^{TOP} (X, Y;\partial)\to
\Cal T^{TOP}(X, \partial X)\to LP_{n-q}(\Delta)
\tag 3.23
$$
and
$$
\cdots \to \Cal S^{TOP} (X, Y;\partial)\to
\Cal S^{TOP}(X, \partial X)\to LS_{n-q}(\Delta)
\tag 3.24
$$

\proclaim{Proposition 1} There exists an $\Omega$@-spectrum
$\Bbb S(X, Y;\partial)$ with homotopy groups
$$
\Cal S_{i}(X, Y;\partial)\cong \pi_i(\Bbb S(X, Y;\partial))\
\text{and}\ \Cal S_{n+1}(X, Y;\partial)\cong \Cal S^{TOP} (X, Y;\partial)).
\tag 3.25
$$
There are algebraic versions
of exact sequences (3.23) and (3.24)
$$
\cdots \to \Cal S_{n+1}(X, Y;\partial)\to
H_n(X, \partial X; \bold L_{\bullet})\overset{\lambda}\to{\to}
 LP_{n-q}(\Delta)\to \cdots
\tag 3.26
$$
 and
$$
\cdots \to \Cal S_{n+1}(X, Y;\partial)\to
\Cal S_{n+1}(X, \partial X)\to LS_{n-q}(\Delta)\to \cdots
\tag 3.27
$$
which are realized on the spectrum level
by cofibrations
$$
(X/\partial X)_+\land\bold L_{\bullet}\to \Sigma^q\Bbb L\Bbb P(\Delta)\to
\Bbb S(X, Y;\partial)
\tag 3.29
$$
and
$$
\Bbb S(X, \partial X)\to \Sigma^{q+1}\Bbb L\Bbb S(\Delta)\to
\Bbb \Sigma S(X, Y;\partial),
\tag 3.30
$$
respectively.
\endproclaim
\demo{Proof} Commutative diagram (2.12)
is generated  by a homotopy commutative diagram of spectra
$$
\CD
@. \vdots @. \vdots @. \vdots \\
@. @VVV @VVV @VVV \\
\cdots @>>>  (\partial X)_+\land \bold L_{\bullet}
@>>>\Bbb L(\rho)
 @>>>
\Bbb S(\partial X) @>>> \cdots
\\
@. @VVV  @VVV  @VVV \\
\cdots @>>> X_+\land \bold L_{\bullet}  @>>> \Bbb L(\pi)
@>>>
\Bbb S^{\partial}(X, \partial X) @>>> \cdots
\\
@. @VVV  @VVV  @VVV\\
\cdots @>>> (X/\partial X)_+\land \bold L_{\bullet} @>>>
\Bbb L(\rho\to \pi) @>>>
\Bbb S(X, \partial X)@>>> \cdots \\
@. @VVV @VVV @VVV \\
@. \vdots @. \vdots @. \vdots
\endCD
\tag 3.31
$$
in which each  row and column is a cofibration  sequence.
Commutative diagram (3.20) is generated
by a homotopy commutative diagram of spectra
$$
\CD
@. \vdots @. \vdots @. \vdots \\
@. @VVV @VVV @VVV \\
\cdots @>>>  \Bbb L\Bbb S(F_{\partial})
@>>>\Bbb L\Bbb P(F_{\partial})
 @>>>
\Sigma^{-q}\Bbb L(\rho) @>>> \cdots
\\
@. @VVV  @VVV  @VVV \\
\cdots @>>> \Bbb L\Bbb S(F)   @>>> \Bbb L\Bbb P(F)
@>>>
\Sigma^{-q}\Bbb L(\pi) @>>> \cdots
\\
@. @VVV  @VVV  @VVV\\
\cdots @>>> \Bbb L\Bbb S(\Delta) @>>>
\Bbb L\Bbb P(\Delta) @>>>
\Sigma^{-q}\Bbb L(\rho\to \pi)@>>> \cdots \\
@. @VVV @VVV @VVV \\
@. \vdots @. \vdots @. \vdots
\endCD
\tag 3.32
$$
in which each  row and column is a cofibration  sequence.
Consider a homotopy commutative square of spectra
$$
\CD
(\partial X)_+\land \bold L_{\bullet} @>>>\Sigma ^q\Bbb L\Bbb P(F_{\partial})\\
@VVV  @VVV   \\
X_+\land \bold L_{\bullet}  @>>> \Sigma ^q\Bbb L\Bbb P(F)\\
\endCD
\tag 3.33
$$
in which the vertical maps are induced by the inclusion $\delta$
and the horizontal maps follows from diagram
(3.10) of
the manifold pair $\partial Y\subset \partial X$ and from
diagram (3.16), respectively.
Denote by $\Bbb S(X,Y; \partial)$ a spectrum fitting
in the diagram extending the square (3.33) by cofibration
sequences
$$
\CD
@. \vdots @. \vdots @. \vdots \\
@. @VVV @VVV @VVV \\
\cdots @>>>  (\partial X)_+\land \bold L_{\bullet}
@>>>\Sigma ^q\Bbb L\Bbb P(F_{\partial})
 @>>>
\Bbb S(\partial X, \partial Y, \partial\xi) @>>> \cdots \\
@. @VVV  @VVV  @VVV \\
\cdots @>>> X_+\land \bold L_{\bullet}  @>>> \Sigma ^q\Bbb L\Bbb P(F)
@>>>
\Bbb S^{\partial}(X, Y, \xi) @>>> \cdots
\\
@. @VVV  @VVV  @VVV\\
\cdots @>>> (X/\partial X)_+\land \bold L_{\bullet} @>>>
\Sigma^q\Bbb L\Bbb P(\Delta) @>>>
\Bbb S(X,Y; \partial)@>>> \cdots \\
@. @VVV @VVV @VVV \\
@. \vdots @. \vdots @. \vdots
\endCD
\tag 3.34
$$
Let $\Cal S_{i}(X, Y;\partial)= \pi_i(\Bbb S(X, Y;\partial))$.
The papers \cite{7}, \cite{11, \S 7.2}, and \cite{14, \S 7.A} 
provide  commutative squares
$$
\matrix  
H_{n+k}(X; \bold L_{\bullet})&\rightarrow
LP_{n-q+k}(\Delta) \\
  t\downarrow\cong & \ \ \downarrow = \\
\Cal T^{TOP}(X\times D^k,\partial (X\times D^k))&
\rightarrow  LP_{n-q+k}(\Delta)\\
\endmatrix
$$
and
$$
\matrix
H_n(X,\partial X; \bold L_{\bullet})&\rightarrow
LP_{n-q}(\Delta) \\
   t\downarrow\cong & \ \ \downarrow = \\
\Cal T^{TOP}(X, \partial X) 
&\rightarrow  LP_{n-q}(\Delta).\\
\endmatrix
$$
Note that the exact sequence (3.26) is obtained by applying the
functor $\pi_0$ to the bottom cofibration in (3.34).
Now, using geometric description of surgery spectra 
(see \cite{7} and \cite{14}), 
an  element $x\in \Cal S_{n+1}(X,Y, \partial)$ 
is defined by a pair 
$(y,z),$  consisting  of a normal map bordism class 
$y\in H_n(X,\partial X; \bold L_{\bullet}),$ for which 
$\lambda(y)=0\in LP_{n-q}(\Delta),$  and  a particular solution 
$z$
of the associated surgery problem for manifold pairs with boundaries
that defines a class of equivalence 
$$
\{(f:M\to X)\}\in \Cal S^{TOP}(X,Y,\partial).
$$
We obtain 
the
map 
$$
\sigma:\Cal S_{n+1}(X,Y, \partial)\to  \Cal S^{TOP}(X,Y,\partial).
$$
Recall, that in geometrically defined exact sequence (3.23) 
the map $LP_{n-q+1}(\Delta)\to \Cal S^{TOP}(X,Y,\partial)$
 is an action.        
Now from the definition of the map $\sigma$ follows 
the commutative diagram (see \cite{7} and \cite{10})
$$
\matrix
LP_{n-q+1}(\Delta)&\to &
\Cal S_{n+1}(X,Y, \partial) &\to & H_{n}(X,\partial X; \bold L_{\bullet})\\
 \downarrow= & &\sigma\downarrow  && t\downarrow \cong\\
LP_{n-q+1}(\Delta) &\to& \Cal S^{TOP}(X,Y,\partial)&\to &
 \Cal T^{TOP}(X, \partial X).
\endmatrix
$$
Using  the Five Lemma we obtain an
isomorphism between (3.23)
and (3.26).
The case of exact sequence (3.27) follows from
a homotopy commutative  diagram of cofibrations
$$
\CD
@. \vdots @. \vdots @. \vdots \\
@. @VVV @VVV @VVV \\
\cdots @>>>
\Bbb S(\partial X, \partial Y, \partial\xi)
@>>>\Bbb S(\partial X)
 @>>>
\Sigma^{q+1}\Bbb L\Bbb S(F_{\partial}) @>>> \cdots \\
@. @VVV  @VVV  @VVV \\
\cdots@>>>
\Bbb S^{\partial}(X, Y, \xi)   @>>> \Bbb S^{\partial} (X, \partial X)
@>>>
\Sigma^{q+1}\Bbb L\Bbb S(F) @>>> \cdots
\\
@. @VVV  @VVV  @VVV\\
\cdots @>>>
\Bbb S(X,  Y; \partial)  @>>>
\Bbb S(X, \partial X)
@>>>
\Sigma^{q+1}\Bbb L\Bbb S(\Delta)@>>> \cdots \\
@. @VVV @VVV @VVV \\
@. \vdots @. \vdots @. \vdots
\endCD
\tag 3.35
$$
which is similar to (3.34). Diagram (3.35) follows from consideration of the
cofibration sequences of the right upper square in (3.35).
\qed
\enddemo
\smallskip

Consider the  composition
$$
LS_*(F_{\partial})\to LP_*(F_{\partial})\to LP_*(F)
\tag 3.36
$$
of geometrically defined maps from diagram (3.20).
The composition (3.36) is realized by a map
of spectra
$$
\Bbb L\Bbb S(F_{\partial})\to \Bbb L\Bbb P(F)
\tag 3.37
$$
which is the composition of maps from diagram (3.32).
We denote the homotopical cofiber of the map (3.37) by
$\Bbb L\Bbb P\Bbb S(\Delta)$
and its homotopy groups by
$$
LPS_i(\Delta)=\pi_i(\Bbb L\Bbb P\Bbb S(\Delta)).
\tag 3.38
$$
In particular, we have a cofibration
$$
\Bbb L\Bbb S(F_{\partial})\to \Bbb L\Bbb P(F)
\to \Bbb L\Bbb P\Bbb S(\Delta).
\tag 3.39
$$
\bigskip

\proclaim{Theorem 2} There exists the following cofibration of spectra
$$
\Bbb S(X, Y, \partial )\to
\Bbb T\Bbb S(X,\partial X)\to \Sigma^{q+1}\Bbb L\Bbb P\Bbb S(\Delta).
\tag 3.40
$$
\endproclaim
\demo{Proof} Consider the following homotopy commutative diagram
of spectra
$$
\CD
\Bbb S(\partial X, \partial Y, \partial \xi)   @>>>
\Bbb S^{\partial}(X, Y, \xi)
@>>>
\Bbb S(X, Y;\partial)\\
 @VVV  @VVV  @VVV\\
\Bbb S(\partial X)   @>>> \Sigma(X_+\land \bold L_{\bullet})
@>>>
\Bbb T\Bbb S(X, \partial X)
\\
 @VVV  @VVV  @VVV\\
\Sigma^{q+1}\Bbb L\Bbb S(F_{\partial})  @>>>
\Sigma^{q+1}\Bbb L \Bbb P(F)
@>>>
\Sigma^{q+1}\Bbb L\Bbb P\Bbb S(\Delta)\\
\endCD
\tag 3.41
$$
in which all rows and columns are cofibration sequences.
The left bottom  square of (3.41) follows from
the commutative diagram
$$
\CD
\Bbb S(\partial X, \partial Y, \partial \xi) @>\simeq>>
\Bbb S(\partial X, \partial Y, \partial \xi) @>>>
\Bbb S^{\partial}(X, Y, \xi)\\
 @VVV  @VVV  @VVV\\
\Bbb S(\partial X)   @>>> \Sigma((\partial X)_+\land \bold L_{\bullet})
@>>>
 \Sigma(X_+\land \bold L_{\bullet}) \\
 @VVV  @VVV  @VVV\\
\Sigma^{q+1}\Bbb L\Bbb S(F_{\partial})  @>>>
\Sigma^{q+1}\Bbb L\Bbb P(F_{\partial})
@>>>
\Sigma^{q+1}\Bbb L \Bbb P(F) \\
\endCD
\tag 3.42
$$
in which the commutative part consisting of the two
 left squares follows from  the diagram 
(3.10) on the spectrum level
for the manifold pair $\partial Y\subset \partial X$ and the commutative part 
consisting of the two 
right squares fits
in (3.34). 
The vertical columns of (3.42) are cofibrations. 
The left column of (3.41) coincides with the left column of (3.42), 
and the middle column of (3.41) coincides with the right column of 
(3.42).
The right vertical maps in (3.41) are defined as map of cofibers of
horizontal maps in accordance with (2.17), (3.35), and (3.39).
The left upper horizontal map in (3.41)
is the composition
$$
\Bbb S(\partial X, \partial Y, \partial \xi)\overset{\simeq}\to{\to}
\Bbb S(\partial X, \partial Y, \partial \xi)
\to
\Bbb S^{\partial}(X, Y, \xi)
\tag 3.43
$$
from the diagram (3.42).
The right  column of cofibration (3.41) is  cofibration (3.40).
\qed\enddemo
\bigskip

\proclaim{Corollary 3} There exists the following long exact sequence
$$
\cdots\to \Cal S_{n}(X,Y;\partial)\to\Cal T\Cal S_n(X, \partial X)\to
LPS_{n-q-1}(\Delta)\to\cdots
\tag 3.44
$$
\endproclaim
\demo{Proof} The exact sequence (3.44) is the homotopy long exact sequence
of the cofibration (3.40).
\qed
\enddemo
\bigskip

\proclaim{Corollary 4} There is a map
$$
\Theta:\Cal T\Cal S^{TOP}(X, \partial X)\to LPS_{n-q}(\Delta)
$$
of obstructions to splitting along the submanifold with boundary
$(Y, \partial Y)$.
For a representative $z=((f,b),(\partial f, \partial b))$ we have
$\Theta (z)=0$ if and only if the class of $z$
in $\Cal T\Cal S^{TOP}(X, \partial X) $ contains a representative
which is split along $(Y, \partial Y)\subset (X, \partial X)$.
\endproclaim
\demo{Proof}
We have the commutative diagram
$$
\CD
\Cal S_{n+1}(X, Y, \partial)   @>>> \Cal S ^{TOP}(X, Y; \partial)\\
@VVV @VVV\\
\Cal T\Cal S_{n+1}(X, \partial X)  @>>>
\Cal T\Cal S^{TOP}(X, \partial X) \\
\endCD
\tag 3.45
$$
in which the horizontal maps are isomorphisms,
the right vertical map is  a natural forgetful map,
and the left vertical map follows from (3.41). Now the diagram
(3.41) and the exact sequence (3.44) provide the result of the Corollary.
\qed
\enddemo
\smallskip

From now on we describe algebraic properties of $LPS_*$@-groups
and their relations to splitting and surgery obstruction groups.
\bigskip

\proclaim{Theorem 5} There exists  a braid of exact sequences
 $$
\matrix \rightarrow &
LP_{n+1}(\Delta) & \rightarrow &
L_{n+q}(\rho) &
\rightarrow &  LS_{n-1}(F_{\partial})& \rightarrow \cr
\ &  \nearrow \ \ \ \ \ \ \ \ \searrow & \ &  \nearrow \ \ \ \ \ \ \
\
\searrow
& \  & \nearrow \ \ \ \ \ \ \ \  \searrow & \ \cr
\ & \ & LP_{n}(F_{\partial})& \ &
LPS_{n}(\Delta) & \ & \ \cr
\ &  \searrow \ \ \ \ \ \ \ \ \nearrow & \ &  \searrow \ \ \ \ \ \ \
\
\nearrow
& \  & \searrow \ \ \ \ \ \ \ \  \nearrow & \ \cr
\rightarrow & LS_{n}(F_{\partial}) & \longrightarrow &
LP_{n}(F) &
\longrightarrow & LP_n(\Delta) & \rightarrow
\endmatrix,
\tag 3.46
$$
 $$
\matrix \rightarrow &
L_{n+q+1}(\pi) & \rightarrow &
LS_{n}(\Delta) &
\rightarrow &  LP_{n}(\Delta)& \rightarrow \cr
\ &  \nearrow \ \ \ \ \ \ \ \ \searrow & \ &  \nearrow \ \ \ \ \ \ \
\
\searrow
& \  & \nearrow \ \ \ \ \ \ \ \  \searrow & \ \cr
\ & \ & L_{n+q+1}(\rho\to \pi)& \ &
LPS_{n}(\Delta) & \ & \ \cr
\ &  \searrow \ \ \ \ \ \ \ \ \nearrow & \ &  \searrow \ \ \ \ \ \ \
\
\nearrow
& \  & \searrow \ \ \ \ \ \ \ \  \nearrow & \ \cr
\rightarrow & LP_{n+1}(\Delta) & \longrightarrow &
L_{n+q}(\rho) &
\longrightarrow & L_{n+q}(\pi) & \rightarrow
\endmatrix,
\tag 3.47
$$
and
 $$
\matrix \rightarrow &
LS_{n}(F_{\partial}) & \rightarrow &
LP_{n}(F) &
\rightarrow &  L_{n+q}(\pi)& \rightarrow \cr
\ &  \nearrow \ \ \ \ \ \ \ \ \searrow & \ &  \nearrow \ \ \ \ \ \ \
\
\searrow
& \  & \nearrow \ \ \ \ \ \ \ \  \searrow & \ \cr
\ & \ & LS_{n}(F)& \ &
LPS_{n}(\Delta) & \ & \ \cr
\ &  \searrow \ \ \ \ \ \ \ \ \nearrow & \ &  \searrow \ \ \ \ \ \ \
\
\nearrow
& \  & \searrow \ \ \ \ \ \ \ \  \nearrow & \ \cr
\rightarrow & L_{n+q+1}(\pi) & \longrightarrow &
LS_{n}(\Delta) &
\longrightarrow & LS_{n-1}(F_{\partial}) & \rightarrow
\endmatrix,
\tag 3.48
$$
where $\rho=\pi_1(\partial X)$ and $\pi=\pi_1(X)$.
Diagrams (3.46) -- (3.48) are realized on the spectrum level.
\endproclaim

\demo{Proof} It follows from the definition of spectra $\Bbb L\Bbb P\Bbb S$
and diagram (3.32) that we have a homotopy commutative diagram of spectra
$$
\matrix
\Bbb L\Bbb S(F_{\partial})& \to &  \Bbb L\Bbb P(F_{\partial})&\to&
\Sigma^{-q}\Bbb L(\rho)\\
\downarrow= &     & \downarrow &&\downarrow\\
\Bbb L\Bbb S(F_{\partial}) & \to    & \Bbb L\Bbb P(F)&\to
&\Bbb L\Bbb P\Bbb S(\Delta)
\endmatrix
\tag 3.49
$$
where the rows are cofibrations and the right vertical
map is induced by the two left vertical maps (see \cite{9}).
The right square in (3.49) is a pullback square since
fibers of horizontal maps are naturally homotopy equivalent.
Hence the right square in
(3.49) is a pushout square and homotopy long exact sequences of
this square give the braid of exact sequences  (3.46).
From diagram (3.32) and \cite{5, Lemma 2} we conclude
that the spectrum
$\Bbb L\Bbb P\Bbb S(\Delta)$ fits in the cofibrations of spectra
$$
\Sigma^{-1}\Bbb L\Bbb P(\Delta)\to \Sigma^{-q}\Bbb L(\rho)\to
\Bbb L\Bbb P\Bbb S(\Delta)
$$
and
$$
\Sigma^{-q-1}\Bbb L(\pi)\to \Bbb L\Bbb S(\Delta)\to
\Bbb L\Bbb P\Bbb S(\Delta).
$$
Now the same line of arguments as for the braid of exact sequences (3.46)
provides diagrams (3.47) and (3.48).
\qed
\enddemo
\bigskip

\proclaim{Corollary 6} The groups $LPD_*(\Delta)$ fit in the following
exact sequences
$$
\cdots \to LS_n(F_{\partial})\to LP_{n}(F)\to LPS_n(\Delta)\to \cdots,
\tag 3.50
$$
$$
\cdots \to LP_{n+1}(\Delta)\to L_{n+q}(\rho)\to LPS_n(\Delta)\to \cdots,
\tag 3.51
$$
and
$$
\cdots \to L_{n+q+1}(\pi)\to LS_{n}(\Delta)\to LPS_n(\Delta)\to \cdots,
\tag 3.52
$$
which are realized on the  level of spectra.
\endproclaim
\demo{Proof} These sequences fit in the diagrams of Theorem 6.
\qed
\enddemo
\bigskip
\proclaim{Corollary 7} Let $\Delta:F_{\partial}\to F$ be an isomorphism
of pushout squares. Then we have isomorphisms
$$
LPS_n(\Delta)\cong L_{n+q}(\rho)\cong L_{n+q}(\pi).
$$
\endproclaim
\demo{Proof} The result follows immediately from the exact sequences
of Corollary 6. \qed
\enddemo

The next theorem describes relations between the obstruction
groups $LPS_*$ and different structure sets which arise naturally
for the manifold pair with boundaries.
\smallskip

\proclaim{Theorem 8}
There exist the following braids of exact sequences
 $$
\matrix \rightarrow &
\Cal  S_{n}(X,Y;\partial) & \rightarrow &
\Cal  T\Cal S_n(X, \partial X) &
\rightarrow &  L_{n-1}(\pi)& \rightarrow \cr
\ &  \nearrow \ \ \ \ \ \ \ \ \searrow & \ &  \nearrow \ \ \ \ \ \ \
\
\searrow
& \  & \nearrow \ \ \ \ \ \ \ \  \searrow & \ \cr
\ & \ & \Cal S_{n}(X, \partial X)& \ &
LPS_{n-q-1}(\Delta) & \ & \ \cr
\ &  \searrow \ \ \ \ \ \ \ \ \nearrow & \ &  \searrow \ \ \ \ \ \ \
\
\nearrow
& \  & \searrow \ \ \ \ \ \ \ \  \nearrow & \ \cr
\rightarrow & L_{n}(\pi) & \longrightarrow &
LS_{n-q-1}(\Delta) &
\longrightarrow & \Cal S_{n-1}(X,Y;\partial) & \rightarrow
\endmatrix,
\tag 3.53
$$
and
  $$
\matrix \rightarrow &
\Cal  S_{n+1}(X,Y;\partial) & \rightarrow &
 H_n(X, \partial X; \bold L_{\bullet}) &
\rightarrow &  L_{n-1}(\rho)& \rightarrow \cr
\ &  \nearrow \ \ \ \ \ \ \ \ \searrow & \ &  \nearrow \ \ \ \ \ \ \
\
\searrow
& \  & \nearrow \ \ \ \ \ \ \ \  \searrow & \ \cr
\ & \ &\Cal  T\Cal S_{n+1}(X, \partial X)& \ &
LP_{n-q}(\Delta) & \ & \ \cr
\ &  \searrow \ \ \ \ \ \ \ \ \nearrow & \ &  \searrow \ \ \ \ \ \ \
\
\nearrow
& \  & \searrow \ \ \ \ \ \ \ \  \nearrow & \ \cr
\rightarrow & L_{n}(\rho) & \longrightarrow &
LPS_{n-q}(\Delta) &
\longrightarrow & \Cal S_{n}(X,Y;\partial) & \rightarrow
\endmatrix,
\tag 3.54
$$
where $\rho=\pi_1(\partial X)$ and $\pi=\pi_1(X)$.
Diagrams (3.53) and (3.54) are realized on the level of spectra.
\endproclaim
\demo{Proof}
Consider a homotopy commutative diagram of spectra
$$
\CD
\Bbb L(\pi)
@>>>\Bbb S(X, \partial X)
 @>>>
\Bbb T\Bbb S(X, \partial X) \\
 @V=VV  @VVV  @VVV \\
\Bbb L(\pi)  @>>> \Sigma^{q+1}\Bbb L\Bbb S(\Delta))
@>>>
\Sigma^{q+1} \Bbb L\Bbb P\Bbb S(X, \partial X) \\
\endCD
\tag 3.55
$$
in which the rows follow from definitions and the columns
are obtained from the natural map of the diagram (3.31) to the diagram (3.32).
The right square of (3.55) is the pullback
and the homotopy long exact sequences of this square
give the commutative diagram (3.53). The case of the diagram (3.54) is similar.
 \qed\enddemo

\bigskip
\subhead 4. Examples
\endsubhead
\smallskip

In this section we compute the $LPS_*$@-groups and natural maps
for several
geometric examples.

Let $(Y^{n-1},\partial Y)\subset (X^n, \partial X), \ n\geq 4$,
be a manifold pair with
boundaries,  where $X$ is a non-trivial $I$@-bundle over the real
projective space $\Bbb RP^{n-1}$ and the submanifold
$Y$ is the restriction of the $I$@-bundle to the
projective space $\Bbb RP^{n-2}\subset \Bbb RP^{n-1}$.
The pair $\partial Y\subset \partial X$ coincides with
$S^{n-2}\subset S^{n-1}$.
We have isomorphisms
$\pi_1(X^n)=\Bbb Z_2^+$ for $n$ odd and
$\pi_1(X^n)=\Bbb Z_2^-$ for $n$ even.
The group $\Bbb Z_2^+$ is a cyclic group of order $2$ 
with the trivial homomorphism
of orientation and $\Bbb Z_2^-$ is this group with 
the nontrivial homomorphism
of orientation.
In the considered case  the squares (1.2) and (3.12)
are the following squares
$$
F^{\pm}=
\left(\matrix
1 &\to& 1\\
\downarrow && \downarrow\\
\Bbb Z_2^{\mp}&\to& \Bbb Z_2^{\pm}\\
\endmatrix
\right)
\tag 4.1
$$
and
$$
F_{\partial}=
\left(\matrix
1\cup 1 &\to& 1\cup 1\\
\downarrow && \downarrow\\
1&\to& 1\\
\endmatrix
\right).
\tag 4.2
$$
All horizontal maps in squares (4.1) and (4.2) are isomorphisms.
\bigskip

\proclaim{Theorem 9} In the considered cases the natural maps
$$
  LP_n(F^{\pm})\to LPS_n(F_{\partial}\to F^{\pm})
$$
fitting in diagrams (3.46) and (3.48) are isomorphisms for
$n=0,1,2,3 \bmod 4$.
Hence we have
$$
LPS_n(F_{\partial}\to F^{+})\cong \Bbb Z_2, \Bbb Z_2, \Bbb Z_2, \Bbb Z
\tag 4.3
$$
and
$$
LPS_n(F_{\partial}\to F^{-})\cong \Bbb Z, \Bbb Z_2, \Bbb Z_2, \Bbb Z_2
\tag 4.4
$$
for $n=0,1,2,3 \bmod 4$, respectively.
\endproclaim
\demo{Proof} It follows from \cite{10, page 153} that
$$
LS_n(F_{\partial})=LN_n(1\cup 1\to 1)=0
$$
for all $n$. From this result and the exact sequence (3.50) the
first statement of the
theorem follows.  We have isomorphisms
$$
LP_n(F^{\pm}) \cong L_{n+1}(i^!_{\mp})
$$
where
$$
i_{\mp}:1\to \Bbb Z_2^{\mp}
$$
is the natural inclusion and
$L_{n+1}(i^!_{\mp})$ is the
relative group of the transfer map
(see, for example, \cite{6}, \cite{8} and \cite{9}).
Now  isomorphisms (4.3) and (4.4) follow
(see,  for example, \cite{9, \S 3} for the case
$F^+$).
\enddemo

\subhead Acknowledgements
\endsubhead
The second author was partially supported by  the Russian
Foundation for Fundamental
Research
grant No. 05--01--00993. The first and third authors were
partially supported by the Ministry of Education, Science and Sport
of the Republic of Slovenia  program No. 509--0101. 
We thank the referee for comments and suggestions.

\vfill

\noindent
Authors' addresses:

\bigskip

\noindent
Matija Cencelj 

\noindent
Institute of Mathematics, Physics and Mechanics, University of
Ljubljana, Jadranska 19, Ljubljana, Slovenia

\noindent
email: matija.cencelj\@fmf.uni-lj.si

\bigskip

\noindent
Yuri V. Muranov

\noindent
Department of Information Science and Management,
Institute of Modern Knowledge,
ulica Gor'kogo 42,
210004 Vitebsk,
Belarus

\noindent
email: ymuranov\@imk.edu.by

\bigskip

\noindent
Du\v san Repov\v s

\noindent
Institute of Mathematics, Physics and Mechanics, University of
Ljubljana, Jadranska 19, Ljubljana, Slovenia

\noindent
email: dusan.repovs\@fmf.uni-lj.si
\enddocument
\bye